\documentclass[twoside]{amsart}

\usepackage{latexsym,amssymb} 

\def\demo{\noindent{\bf Proof. }}
\def\sqr#1#2{{\vcenter{\hrule height.#2pt
        \hbox{\vrule width.#2pt height#1pt \kern#1pt
                \vrule width.#2pt}
        \hrule height.#2pt}}}

\def\build#1_#2^#3{\mathrel{\mathop{\kern 0pt#1}\limits_{#2}^{#3}}}
\def \Box {\hfill\hbox{}\nobreak \vrule width 1.6mm height 1.6mm
depth 0mm \par \goodbreak \smallskip}
\def\rank{{\rm rank}}
\def\th{{\rm th}}
\def\Fitt{{\rm Fitt}}
\def\coker{{\rm coker}}
\def\hgt{{\rm ht}}
\def\dim{{\rm dim}}
\def \Sym{{\rm Sym}}
\def\max{{\rm max}}
\def\Spec{{\rm Spec}}
\def\image{{\rm image}}
\def\p{\varphi}
\def\inc{\subset}
\def\rad{\sqrt}

\newtheorem{Theorem}{Theorem}[section]
\newtheorem{Lemma}[Theorem]{Lemma}

\newtheorem{Proposition}[Theorem]{Proposition}

\begin{document}

\baselineskip=16.5pt

\title[A Simple Proof of Some Generalized Principal Ideal Theorems]
{{\Large\bf A Simple Proof of Some Generalized Principal Ideal Theorems}}

\author[D. Eisenbud, C. Huneke and B. Ulrich]
{David Eisenbud \and Craig Huneke \and Bernd Ulrich}

\thanks{
AMS 1991 {\it Mathematics Subject Classification}. Primary 13C15, 13C40;
Secondary 13D10. \newline\indent
\indent
The authors are grateful to the NSF and to MSRI 
for support.}

\address{MSRI, 1000 Centennial Drive, Berkeley, CA 94720}

\email{de@msri.org}

\address{Department of Mathematics, University of Kansas, Lawrence, KS 66045}

\email{huneke@math.ukans.edu}

\address{Department of Mathematics, Michigan State University, E.
Lansing, MI 48824}

\email{ulrich@math.msu.edu}

\begin{abstract}
Using symmetric algebras we simplify $($and slightly strengthen$)$
the Bruns-Eisenbud-Evans 
``generalized principal ideal theorem'' on the height of order ideals of non 
minimal generators in a module.  We also obtain a simple proof and an 
extension of a result by Kwieci\'nski, which estimates the height of 
certain Fitting ideals of modules having an equidimensional symmetric algebra.
\end{abstract}

\maketitle

\section{Introduction}

The ``generalized principal ideal
theorem'' of Bruns, Eisenbud, and Evans says, roughly, that if
$N$ is a module over a  Noetherian local ring $(R,m)$, and $x$
is an element in $mN$, then the ``order ideal"
$$ N^*(x) = \{\varphi(x) \mid \varphi \in Hom(N,R)\}
$$ 
has height bounded by the rank of $N$. (The usual principal ideal
theorem, that the height of a proper ideal generated by $n$ elements 
is at most
$n$, is the case where $N$ is free of rank $n$.)  The original form of
the result was proved for rings containing a field by Eisenbud-Evans
$($\cite{EE}$)$ and in general by Bruns $($\cite{B}$)$.  
We reduce the statement to the classical principal ideal
theorem, applied in the symmetric algebra of a closely related module.
In this way we obtain a simple proof of a stronger version.

Our proof also yields a new result on the possible height of the ideal
generated by the entries of a row of a matrix whose cokernel $M$ has an
equidimensional symmetric algebra. As a consequence we deduce a result
of Kwieci\'nski $($\cite{K}$)$ which gives a bound for a Fitting ideal
$\Fitt_j(M)$ if
$\sqrt{\Fitt_j(M)}\neq \sqrt{\Fitt_{j+1}(M)}$ 
(Kwieci\'nski's original
proof requires the ring to be a regular affine domain of
characteristic 0).  We also come up with a correction term that allows us to
drop the hypothesis on radicals. 

For further results in the spirit of this paper, see
Johnson $($\cite{J}$)$.

One of the drawbacks of the generalized principal ideal theorem
stated here is the requirement that $x\in mN$. In a subsequent
paper $($\cite{EHU1}$)$ we will show that this requirement can often
be removed or weakened. 

We recall that the applications of the generalized principal ideal
theorem include two useful statements on the {\it determinantal ideal\/}
$I_t(\varphi)$ of $t\times t$ minors of an $n\times m$ matrix
$\varphi$ that generalize the well-known height formulas due to Macaulay
and Eagon-Northcott:

\begin{itemize}
\item[a)] 
If $R$ is local, $I_{t}(\varphi)= 0$ and
$\varphi'$ is the $n\times (m+1)$ matrix derived from
$\varphi$ by adding one column with entries in the maximal ideal, then
$\hgt \, I_{t}(\varphi') \leq n-t+1$ $($Eisenbud-Evans, \cite{EE}$)$.

\item[b)] 
If $I_{t+1}(\varphi) = 0$ then $\hgt \, I_{t}(\varphi) \leq m+n-2t+1$ 
$($Bruns, \cite{B}$)$.
\end{itemize}

\noindent
See the original papers for variations and consequences.

In another paper $($\cite{EHU2}$)$ we will exploit 
symmetric and Rees algebra techniques further to give new bounds
for the height of determinantal ideals under additional
assumptions such as the ring being regular, or the
cokernel of the matrix being torsion free.

For general notions in commutative algebra, see \cite{E}.
A good reference for properties of the symmetric algebra is
\cite{V}.

We started thinking about generalized principal ideal
theorems again because of Kwieci\'nski's interesting
paper \cite{K}, and we would like to thank him for 
sharing it with us.

\section{Symmetric Algebras and Order Ideals}

\begin{Lemma}\label{Lemma1.1}  
Let $R$ be a Noetherian ring and let
$M$ be a finitely generated $R$-module, with free presentation
$$
 F_1 \stackrel{\varphi}{\longrightarrow} F_0 \longrightarrow M \rightarrow 0.
$$
Let $\{T_1, \dots, T_n\}$ be a basis of $F_0$. Let $b:\, F_0 \longrightarrow R$
be any map.  If $P\subset R$ is a prime ideal containing
$\image(b\varphi)$, the ``generalized row ideal of $\varphi$ corresponding to $b$,''
then
$$ \dim \; R_P/\image(b\varphi) \geq
 \dim \; \Sym(M)_{(P,\; \{T_i-b(T_i)\})} -n .
$$
\end{Lemma}
\demo The equality
$$ R_P/(\image(b\varphi)) = \Sym(M)_{(P,\;
\{T_i-b(T_i)\})}/(\{T_i-b(T_i)\})
$$
shows that the $n$-generator ideal 
$(\{T_i-b(T_i)\})$ is proper.
The dimension estimate follows since the ring
$\Sym(M)_{(P,\; \{T_i-b(T_i)\})}$ is local $($see
\cite[Corollary~10.9]{E}$)$.
\Box

\medskip

For an $R$-module $M$ and a prime $P\subset R$ we write $\mu_P(M)$ for
the minimal number of generators of the
$R_P$ -module $M_P$. If $(R,m)$ is local we set $\mu(M)=\mu_m(M)$.

Here is the promised strengthening of the generalized principal
ideal theorem of Bruns, Eisenbud, and Evans $($\cite[Theorem~1]{B}$)$:

\begin{Theorem}\label{Theorem1.2}  
Let $(R,m)$ be a Noetherian local ring and let
$N$ be a finitely generated $R$-module. If $x\in mN$, then
$\dim \, R/N^*(x) \geq \dim \, R/Q-\mu_Q(N)$ for every minimal prime $Q$
of
$R$; thus for every such $Q$,
$$ 
\hgt \, N^*(x) \leq \dim \, R-\dim \, R/Q +\mu_Q(N).
$$
\end{Theorem}

Note that if $R$ is a domain then
the right hand side of the inequality is the rank of $N$.
The improvement over the original versions is in the treatment of
the rank in the general case.
For the proof we need:

\begin{Lemma}\label{Lemma1.3}  
Let $R$ be a  Noetherian local ring, let $N$ an
$R$-module, and let $\pi:\, F\longrightarrow N$ be a surjection from a finitely
generated free module. If $M = \coker(\pi^*)$, then
$$
\mu(F) \leq \mu(M)+\mu(N).
$$ 
\end{Lemma}

\noindent{\bf Proof of Lemma~\ref{Lemma1.3}.\/} 
The kernel of $\pi$ contains a free
summand $G$ of $F$ with rank $\mu(F)-\mu(N)$. The projection map
$F^*\longrightarrow G^*$ factors through $M$, so 
$\mu(M)\geq \mu(G^*) = \mu(F)-\mu(N)$. \Box

\medskip

\noindent{\bf Proof of Theorem~\ref{Theorem1.2}.\/}  Let
$$
F \stackrel{\pi}{\longrightarrow} N \rightarrow 0
$$ 
be a surjection from a free module
$F$ of rank $n = \mu(N)$ onto $N$, and let $M=\coker(\pi^*)$.  Let
$b\in mF$ be an element mapping onto $x$; we regard $b$ as a map from
$F^*$ to $R$. Choose a free presentation of $M$ of the form
$$ 
F_1\stackrel{\pi}{\longrightarrow} F^*\longrightarrow  M \rightarrow 0.
$$ 
The order ideal of $x$ is given by
$N^*(x) = \image (b\varphi)$, so we may apply Lemma~\ref{Lemma1.1}.  Since all
$b(T_i)\in m$ we get
$$ 
\dim \, R/N^*(x) \geq \dim \, \Sym(M)_{(m,\{T_i\})}-n.
$$ 
As $\Sym(M)$ is graded we may use the dimension formula of Huneke
and Rossi $($\cite[Theorem~2.6]{HR}, \cite[Theorem~1.2.1]{V}$)$ 
to obtain
$$ 
\dim \, \Sym(M)_{(m,\{T_i\})} = \dim \, \Sym(M) = \max\{
\dim \, R/Q +\mu_Q(M)
\mid Q\in \Spec(R) \}.
$$ 
By Lemma 1.3, $\mu_Q(M) \geq n-\mu_Q(N)$. Combining these formulas, we
conclude that
$\dim \, R/N^*(x) \geq \dim \, R/Q-\mu_Q(N)$ for all primes $Q$, and
hence in particular for all minimal primes, as required for the first
statement.  The second statement  follows at once.\Box

\section{Heights of Determinantal Ideals}

The formulas for the heights of determinantal ideals that follow from
Theorem~\ref{Theorem1.2} are often sharp. But under an additional hypothesis 
on $M$, Kwieci\'nski $($\cite[Theorem~1]{K}$)$ 
gave a remarkable new bound on the
heights of the Fitting ideals of $M$. (Recall that with notation
as in Lemma~\ref{Lemma1.1}, the $i^\th$ {\it Fitting ideal\/}
$\Fitt_i(M)$ is  $I_{n-i}(\varphi)$.) His proof required $R$ to
be a regular affine domain over a field of characteristic 0, 
but using Lemma~\ref{Lemma1.1} we can remove this restriction and prove
a stronger result.

We deduce our theorem from a more precise version, giving
an estimate for the heights of the generalized row ideals of $\varphi$:

\begin{Theorem}\label{Theorem2.1} 
Let $R$ be a Noetherian universally catenary, locally equidimensional ring
whose maximal
ideals all have the same height. Let $M$ be a finitely generated
$R$-module such that
$\Sym(M)$ is equidimensional.  If $P$ is a minimal prime of the
generalized row ideal $\image(b\varphi)$   of Lemma 1.1, then
$\hgt \, P \leq n+\hgt\, Q -\mu_Q(M)$ for any prime $Q\subset P$.
\end{Theorem}

Informally one may state this result as: under the given
hypothesis, the row ideals behave as though $M$ had 
projective dimension 1. For the proof we need the next proposition.

\begin{Proposition}\label{Proposition2.2}  
Let $R$ be a Noetherian universally catenary, locally equidimensional 
ring whose maximal ideals all have the same height.
If $A$ is a  graded, equidimensional $R$-algebra, finitely
generated by elements of positive degrees, and $P\subset R$ is a prime
ideal, then the maximal ideals of
$A_P$ that contract to $P$ all have height equal to $\dim \, A_P$.
Moreover, $A_P$ is equidimensional.
\end{Proposition}
\demo
We may write $A = S/J,$ where
$S$ is a positively graded polynomial ring $R[T_1,\dots,T_n]$ and 
$J \subset S$ is a homogeneous ideal.

We shall first show that $J$ is height-unmixed.  Let $Q$ be a
(necessarily homogeneous) minimal prime of $J$. Since $S$ is positively 
graded, $J$ is homogeneous, and $A$ is 
equidimensional, there exists a homogeneous maximal  ideal ${\mathcal M}$  of $S$
such that
$\hgt \, {\mathcal M}/Q = \dim \, A$.   As ${\mathcal M}$ is homogeneous, it has the form
${\mathcal M} = mS+S_+$ where $m$ is a maximal ideal of $R$. Hence by our
assumptions on $R$ we have $\hgt \, {\mathcal M} = \dim \, R+n$. Since 
$S_{\mathcal M}$ is
equidimensional and catenary,  $ \hgt \, Q = \hgt \, {\mathcal M} - 
\hgt \, {\mathcal M}/Q = \dim \, R + n -\dim \, A$.

As $J$ is height-unmixed, $J_P\subset S_P$ has the same property.
Therefore  $A_P = S_P / J_P$ is equidimensional, because $S_P$ is
a catenary, equidimensional positively graded ring over a local ring and 
$J_P$ is a homogeneous ideal.

Next let ${\mathcal M} \subset S$ be  maximal among the ideals of  $S$ 
that contract to $P$.  Since $S_{\mathcal M}$ is equidimensional of 
dimension $\dim \, R_P + n$ and $J$ is height-unmixed, we obtain 
$\dim \, A_{\mathcal M} = \dim \, S_{\mathcal M} - \hgt \, J_{\mathcal M} = 
\dim \,  R_P + n - \hgt \, J$.  Hence all maximal ideals of $A_P$ 
contracting to $P$ have the same height.  On the other hand, taking 
${\mathcal M} = PS + S_+$, we have $\dim \,
A_{\mathcal M} = \dim \, A_P$, because 
$A_P$ is positively graded over the local ring $R_P$.
\Box
\bigskip

\noindent{\bf Proof of Theorem 3.1.\;}
By Proposition \ref{Proposition2.2}, 
$$
\dim\ \Sym(M)_{(P, \{T_i-b(T_i)\})} =\dim\ \Sym(M)_P.
$$
On the other hand, the Huneke-Rossi dimension formula
$($\cite[Theorem 2.6]{HR}$)$ implies that
$$
\dim\ \Sym(M)_P \geq \dim(R/Q)_P+\mu_Q(M).
$$
Now we can apply Lemma \ref{Lemma1.1}.
\Box

\bigskip

Here is our generalization of Kwieci\'nski's Theorem:

\begin{Theorem}\label{Theorem2.3} 
Let $R$ be a regular domain whose maximal ideals
all have the same height. Let $M$ be a finitely generated
$R$-module such that $\Sym(M)$ is equidimensional.  Let $i$ be a
positive integer. If $P$ is a minimal prime ideal of the Fitting
ideal
$\Fitt_{i-1}(M)$ not containing
$\Fitt_i(M)$, then 
$$\hgt \, P \leq i(i-\rank(M)).$$
\end{Theorem}
\demo
By Proposition 2.2 the algebra $\Sym(M)_P$ is
again equidimensional, so we may replace $R$ by
$R_P$, and thus assume that $i = \mu_P(M)$. Taking a free presentation
of
$M$ as in Lemma 1.1, we must show that the ideal $I_1(\varphi)$ generated
by the entries of a matrix representing $\varphi$ has height $\leq
i(i-\rank(M))$.  This ideal is the sum of the $i$ ideals generated by
the entries of the rows of a matrix representing $\varphi$. By Theorem 2.1
each of these ideals has height $\leq i-\rank(M)$.  By Serre's Theorem
$($\cite[Theorem 3, p.~V--18]{S}$)$ the heights of ideals in regular
local rings are subadditive, and the desired inequality follows.
\Box

\bigskip

The height estimate of Theorem~\ref{Theorem2.3} does not extend to the 
case in which
$\rad{{\rm Fitt}_{i-1}(M)} = \rad{{\rm Fitt}_{i}(M)}$. For instance, if $M$ 
is a complete intersection ideal in a regular local ring with $n =
\mu(M)\geq 3$ and  $i = 2$, then $\hgt\, {\rm Fitt}_{i-1}(M) = n > 2 =
i(i-\rank(M))$.
Our next estimate explains and covers this situation. 
It refines Theorem~\ref{Theorem2.3} by adding a correction term
for the case when the two Fitting ideals have the same radical.

\begin{Theorem}\label{Theorem2.4}  
Let $R$ be a regular domain whose maximal ideals
all have the same height. Let $M$ be a finitely generated
$R$-module such that $\Sym(M)$ is equidimensional.  Let $i\geq \rank(M)$ 
be an
integer. If $P$ is a minimal prime ideal of the Fitting
ideal
$\Fitt_{i-1}(M)$, then 
$$\hgt \, P \leq i(i-\rank(M))+\mu_P(M)-i.$$
\end{Theorem}
\demo
By Proposition~\ref{Proposition2.2} we may again replace $R$ by $R_P$. For 
$e=\rank(M),\ n =
\mu(M)$ and $\p$ a matrix with $n$ rows presenting $M$, we are going to prove
that $\hgt\, I_{n-i+1}(\p) \leq i(i-e)+n-i$.

We may assume that
$I_{n-i+1}(\p)\ne
0$.
For $n-i\leq j\leq n$ let $\p_j$ be the matrix consisting of the first $j$
rows
of $\p$ and define
$$ t = {\rm min}\{j|\, I_{n-i+1}(\p)\inc \rad{I_{n-i+1}(\p_j)}\,\}.$$
Obviously $n-i+1\leq t\leq n$.

By Theorem~\ref{Theorem2.1} the $t^\th$ row ideal $J_t$ of $\p$ has height at most
$n-e$ and hence by the definition of $t$,
$$
\hgt\, I_{n-i+1}(\p) \leq n-e + \hgt\, I_{n-i+1}(\p_{t-1}). 
\leqno{(*)}
$$
There exists a prime ideal $Q$ with $I_{n-i+1}(\p_{t-1})\inc Q$ and
$I_{n-i+1}(\p_t)\not\subset Q$. Since 
$I_{n-i+1}(\varphi_t)\subset I_{n-i}(\varphi_{t-1})$ we automatically
have $I_{n-i}(\varphi_{t-1})\not\subset Q$.
Without changing $I_{n-i+1}(\p_{t-1})_Q$ we
may therefore assume
that over the ring $R_Q$,
$$
\varphi_Q=
\left[
\vcenter{
\hbox{%
\hbox to .5in{\hss$1_{n-i}$\hss}%
\vrule height .3in depth .2in
\hbox to 1.0055in{\hss$0$\hss}%
}
\hrule
\hbox{$
\vcenter{\hbox to .5in{\hss$0$\hss}}
\vrule
\vcenter{
\hbox{%
\hbox to .5in{\hss$\varphi'$\hss}%
\vrule height .3in depth .2in
\hbox to .5in{\hss$\varphi''$\hss}%
}
\hrule
\hbox{%
\hbox to .5in{\hss$0$\hss}%
\vrule height .3in depth .2in
\hbox to .5in{\hss$1_{s-n+i}$\hss}%
}}$}}
\right]
$$
where $n-i+1\leq s\leq n$ and $\p', \p''$ have $n-s$ rows and entries in the
maximal ideal of $R_Q$.

The symmetric algebra of Coker$(\p')\cong M_Q$ is equidimensional by
Proposition~\ref{Proposition2.2}, and therefore Theorem~{Theorem2.3} 
(or \ref{Theorem2.1}) yields
$\hgt\, I_1(\p')\leq (n-s)(n-s-e)$. On the other hand, $\hgt\, I_1(\p'')\leq
\mu(I_1(\p''))\leq (n-s)(s-n+i)$. As $I_{n-i+1}(\p_{t-1})_Q\inc I_1(\p') +
I_1(\p'')$, inequality $(*)$ implies $\hgt\, I_{n-i+1}(\p)\leq n-e +
(n-s)(n-s-e) + (n-s)(s-n+i)
\leq i(i-e)+n-i$.
\Box

\end{document}